\documentclass{article}
\usepackage{color}

\usepackage{ifthen}
\usepackage[latin1]{inputenc}
\usepackage{amsmath,amssymb,amsthm}
\usepackage[dvips]{graphicx,epsfig}
\usepackage{color} %pour xfig
\newcommand{\A}{\mathbb{A}}
\def\AA{{\mathbb A}}

\def\GG{{\mathbb G}}

\def\NN{{\mathbb N}}

\def\PP{{\mathbb P}}

\def\fxb{{\ensuremath \mathcal B}}
\def\fxc{{\ensuremath \mathcal C}}

\def\fxl{{\ensuremath \mathcal L}}

\def\fxo{{\ensuremath \mathcal O}}

\newcommand{\fd}{\ensuremath{\rightarrow}}

\newcommand{\findem}{\nolinebreak\vspace{\baselineskip} \hfill\rule{2mm}{2mm}\\}
\renewcommand{\phi}{\ensuremath{\varphi}}
\newcommand{\inc}{\ensuremath{\subset}}

\newcommand{\ox}{\otimes }
\newcommand{\pla}{\A^2}
\newcommand{\plp}{\PP^2}
\renewcommand{\phi}{\ensuremath{\varphi}}
\newcommand{\x}{\ensuremath{\times}}

\newcommand{\s}{Spec\ }

\newtheorem{nt}{Notation}

\newtheorem{prop}[nt]{Proposition}

\newtheorem{exercice}[nt]{Exercice} 
\newcounter{numeroquestion} 
\newcounter{numerosousquestion}

\newcommand{\sousquestion}{\ifthenelse{\value{numerosousquestion}=1}{}{\\}\textbf{\roman{numerosousquestion})} \addtocounter{numerosousquestion}{1}}

\newenvironment{listecompacte}
{\begin{list}
    {\ensuremath{\bullet}}
    {\setlength{\topsep}{2pt}
      \setlength{\itemsep}{1pt} \setlength{\parsep}{0pt}}
}
{\end{list}
}

\newtheorem{coro}[nt]{Corollary} 
\newtheorem{defi}[nt]{Definition}

\newtheorem{lm}[nt]{Lemma} 
 
\newtheorem{rem}[nt]{Remark} 
\newtheorem{thm}[nt]{Theorem}  
 
\newcommand{\demo}{\noindent \textit{Proof}}

\begin{document}
\sloppy
\title{Computing limit linear series with infinitesimal methods}
\date{}
\author{Laurent Evain (laurent.evain@univ-angers.fr)}
\maketitle

%%%%%%%%%%%%%%%%%%%%%%%%%%%%%%%%%%%%%%%%%%%%%%%%%%%%%%%%%%%%%%%%
%%%%                    Debut du contenu                  %%%%%
%%%%%%%%%%%%%%%%%%%%%%%%%%%%%%%%%%%%%%%%%%%%%%%%%%%%%%%%%%%%%%%%

\section*{Abstract: }
Alexander and Hirschowitz 
\cite{alexander-hirschowitz00:grosPointsSuffisammentNombreux}
determined the Hilbert function of a generic
union of fat points in a projective space 
when the number of fat points is much bigger than 
the greatest multiplicity of the fat points. Their method is based
on a lemma which determines the limit of a linear system depending on fat
points which approach a divisor. \\
On the other hand, Nagata \cite{nagata59:14emeProblemeDeHilbert}, 
in connection with its counter
example to the fourteenth problem of Hilbert determined the 
Hilbert function $H(d)$ of the union of $k^2$ points of the same multiplicity
$m$ in the plane up to degree $d=km$. \\
We introduce a new method to determine limits of linear
systems.
This generalizes the result by Alexander and Hirschowitz. 
Our main application of this method is the conclusion of the work
initiated by Nagata: we
compute $H(d)$ for all $d$.
As a second application, we determine the generic successive collision 
of four fat points of the same multiplicity in the plane.

\section{Introduction}
\label{sec:introduction}
Let $X$ be a (quasi-)projective scheme,  $\fxl$ a linear system on $X$ and 
$Z\inc X$ a generic 0-dimensional subscheme. 
In this paper, we adress the problem of 
determining the dimension of $\fxl(-Z)$, or more precisely the 
limit of $\fxl(-Z)$ when $Z$ specializes to a subscheme $Z'$.
\\
Our result gives an estimate of this limit
when $Z$ moves to a divisor and 
satisfies suitable conditions( $Z$ is 
the generic embedding of a union 
$Z_1 \cup Z_2 \dots \cup Z_s$ of  
monomial schemes).
More precisely, we introduce 
a combinatorical procedure to construct a 
system $\fxl'$, ``simpler'' than $\fxl$ in the sense 
that it has smaller degree, and we settle an inclusion 
$\lim \fxl(-Z)\subset \fxl'$. In concrete exemples (see 
the applications below), the inclusion suffices to 
compute $\dim \fxl(-Z)$: there is an expected dimension 
$d_e$ which verifies 
\begin{displaymath}
d_e\leq \dim \fxl(-Z)= \dim \lim \fxl(-Z)
\leq \dim \fxl'=d_e,
\end{displaymath}
hence $\dim \fxl(-Z)=d_e$.
\\
To give a flavour of the theorem, suppose for
simplicity that $Z$ is the generic fiber 
of a subscheme $F \inc X\x \AA^1$ 
flat over $\AA^1=\s k[t]$ and such 
that the support of the fiber $F(t)$ 
approaches a divisor $D$ when $t\fd 0$. We find an 
integer $r$ and a residual scheme $F_{res}\inc F(0)$
such that 
\begin{displaymath}
\lim_{t\fd 0}\fxl(-F(t)) \inc \fxl(-rD-Z_{res}).
\end{displaymath}
There is a trivial inclusion
\begin{displaymath}
\lim_{t\fd 0}\fxl(-F(t)) \subset \fxl(-F(0)),
\end{displaymath}
but of course our result is more detailed and is not reductible to
this trivial case. In the examples we consider, the last
inclusion of the tower
\begin{displaymath}
\lim_{t\fd 0}\fxl(-F(t))\subset  \fxl(-rD-Z_{res}(0)) \inc \fxl(-F(0))
\end{displaymath}
is always a strict inclusion.
\\
The method to prove the result is infinitesimal in nature.
There is a unique flat family $G$ over 
$\AA^1$ whose fiber over a general $t\neq 0$ is $\fxl(-F(t))$.  
Our theorem is obtained with a careful analysis of the restrictions 
$G\x_{\AA^1}\s k[t]/(t^{n_i}) \inc G$ for well chosen 
integers  $n_1,\dots,n_r$. 
\\
Our theorem generalizes the main lemma of Alexander-Hirschowitz
\cite{alexander-hirschowitz00:grosPointsSuffisammentNombreux}.
Their statement corresponds essentially 
to ours in the special case $r=1$. However, the proofs are 
different. In fact, when Alexander-Hirschowitz published their 
theorem, our theorem did already exist in a weaker version where 
the $0$-dimensional subscheme $Z$ moving to the divisor had to
be supported by a unique point. The current version is a merge 
which contains both our earlier version and Alexander-Hirschowitz 
version. 
\medskip \\
As an application of our theorem, we extend results by Nagata
relative to the Hilbert functions of fat points in the plane. 
In connection with his construction of the counter example to 
the fourteenth problem of Hilbert, Nagata proved that 
the Hilbert function of a generic union $Z$ of $k^2$ fat points
of the same multiplicity $m$ in $\plp$ is $H_Z(d)=\frac{(d+1)(d+2)}{2}
$ if the degree is not to big, namely if $d\leq km$. This result 
is asymptotically optimal in $m$ in the sense that it is sufficient 
to compute the Hilbert function up to the critical 
degree $d=km+[ \frac{k}{2} ]$ to determine the whole Hilbert function.
Nagata was just missing the last extreme hardest $[\frac{k}{2}]$
cases. 
We compute 
the Hilbert function for every degree: 
$H_Z(d)=min(\frac{(d+1)(d+2)}{2},k^2 \frac{m(m+1)}{2})$.
This result was already proved when the number of points is a power 
of four in \cite{evain99:4hgrosPoints} by methods relying on 
the geometry of integrally closed ideals which we could not 
push further. 
\\
Putting the result in 
perspective, we recall that a consequence of  Alexander-Hirschowitz
\cite{alexander-hirschowitz00:grosPointsSuffisammentNombreux}
is that the Hilbert function of a generic union of $k$ fat points 
in the plane of multiplicity $m_1,\dots,m_k$ is 
$H_Z(d)=min(\frac{(d+1)(d+2)}{2},\sum_{i=1}^{k} \frac{m_i(m_i+1)}{2})$
provided $k>>max(m_i)$. In view of their result, we are left with the 
cases when the multiplicities are not too small with respect to the 
number of points. Among these, it is known empirically that the 
hardest cases are those  with a fixed number of points and big
multiplicities. Our theorem includes such cases. 
\\
As a second application, we compute the generic successive collision
of four fat points in the plane of the same multiplicity (recall that 
a \textit{successive} collision of punctual
schemes  $Z_1, \dots,Z_s$ is a subscheme obtained as a flat limit
when the $Z_i$'s approach one after the other, ie. you
first collide $Z_1$ and $Z_2$ in a subscheme $Z_{12}$, then you collide $Z_3$ with the
previous collision $Z_{12}$ and so on... A \textit{generic} successive
collision is a successive collision where by definition the $Z_i$'s
move on generic curves of high degree ). \\
Let us explain the motivations for such a computation. 
First, collisions determine the Hilbert function 
of the  generic union $Z$ of the fat points. Indeed, 
there exist ``universal'' collisions $C_0$ on which one can read off 
the Hilbert function of $Z$:
$\forall d,\ H_Z(d)=H_{C_0}(d)$
\cite{evain97:lesCollisionsDeterminentLaPostulation-CRAS}.
Moreover, constructing collisions is a useful technical
tool of the Horace method 
(see \cite{hirscho85:methodeHoraceManuscripta}).
\\
However, determining all collisions of any number of fat points 
is far beyond our knowledge 
since this problem is far more difficult than the open and long standed 
problem of determining the Hilbert function of a generic union
of fat points. It is thus natural to restrict our attention to 
special collisions. In view of the postulation problem, 
one looks for collisions special enough so that it is possible to
compute them, but general enough so that they can stand for a
universal collision in the above sense. A natural class of collisions
to be considered is the class of generic successive collisions. 
Can we compute them ? Is there a universal collision among them ?  
A generic successive collision of three fat points is universal
\cite{evain_these}, ie. this collision has the same Hilbert function 
as the generic union of the three fat points. 
We use our theorem to compute the generic successive collision of 
four fat points. 
Our computation proves that this collision is not universal.
Beyond this example, the computation also illustrates how our theorem 
can be used to determine many collisions, thus 
extending the toolbox of the Horace method.

\section{Statement of the theorem}
\label{sec:statement-theorem}
We fix a generically smooth quasi-projective scheme $X$ of dimension
$d$, a locally free sheaf $L$ of rank one on $X$ and a sub-vector space
$\fxl\inc H^0(X,L)$. Let $Z\inc X_{k(Z)}$ be a $0$-dimensional 
subscheme parametrised by a 
non closed point of  $Hilb(X)$ with residual field $k(Z)$.
Let  $\fxl(-Z)\inc \fxl$ be the sub-vector space of sections 
which vanish on $Z$ (see the definition below). 
Our goal is to give an estimate of the dimension
$\dim \fxl(-Z)$
under suitable conditions. 

A staircase $E\subset\NN^d$ is a subset whose complement 
$C=\NN^d\setminus E$ verifies $\NN^d + C \subset C$. We denote 
by $I^E$ the ideal of $k[x_1,\dots,x_d]$ (resp. of
$k[[x_1,\dots,x_d]])$, of $k[[x_1,\dots,x_d]][t]\dots$)
generated by the monomials $x_1^{e_1}.\dots .x_d^{e_d}=x^e$ 
whose exponent $e=(e_1,\dots,e_d)$ is 
in $C$. If $E$ is a finite staircase, the subscheme $Z(E)$ defined by $I^E$
is 0-dimensional and its degree is $\# E$. The
map $E\mapsto Z(E)$ is a one-to-one correspondance between the finite 
staircases of $\NN^d$ and the monomial punctual subschemes of $\s
k[x_1,\dots,x_d]$. If $E=(E_1,\dots,E_s)$ is a set of finite 
staircases, if $X$ is irreducible 
and if $Z(E)$ is the (abstract non embedded) 
disjoint union $Z(E_1)\coprod \dots \coprod Z(E_s)$,
there is an irreducible scheme $P(E)$ which parametrizes the 
embeddings $Z(E)\fd X_s$, where $X_s\inc X$ is the smooth locus
(\cite{grothendieck60:techniquesConstructionEtSchemasDeHilbert} and 
\cite{hirscho85:methodeHoraceManuscripta}). 
Such an embedding  $Z(E)\fd X_s$ determines a subscheme of $X$, thus 
there is a natural morphism $f:P(E)\fd Hilb(X)$ to the Hilbert scheme of
$X$. We denote by $X(E)$ the subscheme 
parametrised by $f(p)$ where $p$ is 
the generic point of $P(E)$. We will say that $X(E)$ is the 
generic union of the schemes $Z(E_1),\dots,Z(E_n)$.  
If $Z\inc X$ is a subscheme, denote by 
$\fxl(-Z) \inc \fxl$ the subvector space which contains the elements 
of $\fxl$ vanishing on $Z$. If $p$ is a non closed point of
$Hilb(X)$ whose residual field is $k(p)$, and if $Z \inc X\x_k \s k(p)$ is
the corresponding subscheme, the definition of $\fxl(-Z)$ is as follows. 
Since $\fxl \ox k(p) \inc H^0(L\ox k(p),X\x k(p))$, it makes sense
to consider the vector space $V \inc \fxl \ox k(p)$ containing the
sections which vanish on $Z$. Denoting by $\lambda$ the codimension of $V$,
we may associate with $V$ a $k(p)$-point $g \in Grass_{k(p)}
(\lambda, \fxl \ox k(p))= Grass_k(\lambda,
\fxl) \x \s k(p)$ (\cite{grothendieck60:ega1}, prop.9.7.6). 
In particular $\fxl(-Z)$ is well defined as a (non closed) point of 
$Grass_k(\lambda,\fxl)$. The goal of the theorem is to give 
an estimate of $\dim \fxl(-X(E))$. \\
To formulate the theorem, we need some combinatorial notations that we 
introduce now.  The $k^{th}$ slice of a staircase 
$E\inc \NN^d$ is the staircase $T(E,k)\inc \NN^d$ 
defined by: 
\begin{displaymath}
  T(E,k)=\{(0,a_2,\dots,a_d)\ such \ that\ (k,a_2,\dots,a_d)\in E\}
\end{displaymath}
If $E=(E_1,\dots,E_s)$ is a s-tuple of staircases and
$t=(t_1,\dots,t_s)$,
we set 
\begin{displaymath}
  T(E,t)=(T(E_1,t_1),T(E_2,t_2),\dots,T(E_s,t_s)).
\end{displaymath}
A staircase $E\inc \NN^d$ is characterized by a height function
$h_E:\NN^{d-1}\fd \NN$ which verifies: 
\begin{displaymath}
\forall a,b \in \NN^{d-1}, \ h_E(a+b)\leq h_E(a)
\end{displaymath}
The staircase $E$ and $h_E$ can be deduced one from the other via the 
relation:
\begin{displaymath}
  (a_1,\dots,a_d)\in E \Leftrightarrow a_1<h_E(a_2,\dots,a_n)
\end{displaymath}
The staircase $S(E,t)$ is defined by its height function:
\begin{eqnarray*}
h_{S(E,t)}(a_2,\dots,a_d)&=&h_E(a_2,\dots,a_d)\mbox{ if }
t \geq  h_E(a_2,\dots,a_d)
\\
&=&h_E(a_2,\dots,a_d)-1 \mbox{ if } t < h_E(a_2,\dots,a_d)
\end{eqnarray*}
Intuitivly, it is the staircase obtained by the suppression of 
the $t^{th}$ slice, as shown by the following figure.
\\
\begin{center}

\setlength{\unitlength}{0.00083300in}%
\begingroup\makeatletter\ifx\SetFigFont\undefined%
\gdef\SetFigFont#1#2#3#4#5{%
  \reset@font\fontsize{#1}{#2pt}%
  \fontfamily{#3}\fontseries{#4}\fontshape{#5}%
  \selectfont}%
\fi\endgroup%
\begin{picture}(2274,1856)(6889,-1605)
\thicklines
\put(6901,-699){\line( 1,-1){187.500}}
\put(7089,-886){\line( 0,-1){188}}
\put(7464,-1074){\line( 0, 1){188}}
\put(7464,-886){\line( 1, 1){187}}
\put(7464,-886){\line( 0, 1){187}}
\put(7464,-699){\line( 1, 1){187.500}}
\put(7089,-886){\line( 1, 1){187}}
\put(7276,-699){\line( 0,-1){187}}
\put(7089,-511){\line( 1,-1){375}}
\put(7276,-699){\line( 0, 1){188}}
\put(7276,-511){\line( 1, 1){187.500}}
\put(7089,-324){\line( 1,-1){375}}
\put(7089,-136){\line( 1,-1){187.500}}
\put(7276,-324){\line( 1, 1){188}}
\put(7089, 51){\line( 1,-1){187}}
\put(7276,-136){\line( 1, 1){187.500}}
\put(7276,-136){\line( 0,-1){375}}
\put(7089,-324){\line( 0,-1){187}}
\put(7089,-511){\line( 1,-1){375}}
\put(7464,-886){\line( 0, 1){187}}
\put(7464,-699){\line(-1, 1){375}}
\put(7464,-699){\line( 1, 1){187.500}}
\put(7651,-511){\line( 0,-1){188}}
\put(7651,-699){\line(-1,-1){187}}
\put(7276,-511){\line( 1, 1){187.500}}
\put(7464,-324){\line( 1,-1){187}}
\put(7651,-511){\line(-1,-1){187.500}}
\put(7464,-699){\line(-1, 1){188}}
\put(6901,-699){\line( 1, 1){188}}
\put(7089,-511){\line( 0, 1){562}}
\put(7089, 51){\line( 1, 1){187.500}}
\put(7276,239){\line( 1,-1){188}}
\put(7464, 51){\line( 0,-1){375}}
\put(7464,-324){\line( 1,-1){187}}
\put(7651,-511){\line( 0,-1){375}}
\put(7651,-886){\line(-1,-1){187.500}}
\put(7464,-1074){\line(-1, 1){188}}
\put(7276,-886){\line(-1,-1){187.500}}
\put(7089,-1074){\line(-1, 1){188}}
\put(6901,-886){\line( 0, 1){187}}
\put(7276,-511){\line( 0,-1){188}}
\put(7464,-699){\line( 0,-1){187}}
\put(8401,-699){\line( 1, 1){188}}
\put(8589,-511){\line( 0, 1){375}}
\put(8589,-136){\line( 1, 1){187}}
\put(8776, 51){\line( 1,-1){187.500}}
\put(8964,-136){\line( 0,-1){375}}
\put(8964,-511){\line( 0, 1){  0}}
\put(8964,-511){\line( 1,-1){187.500}}
\put(9151,-699){\line( 0,-1){187}}
\put(9151,-886){\line(-1,-1){187.500}}
\put(8964,-1074){\line(-1, 1){188}}
\put(8776,-886){\line(-1,-1){187.500}}
\put(8589,-1074){\line(-1, 1){188}}
\put(8401,-886){\line( 0, 1){187}}
\put(8401,-699){\line( 1,-1){187.500}}
\put(8589,-886){\line( 0,-1){188}}
\put(8589,-886){\line( 1, 1){187}}
\put(8776,-699){\line( 0,-1){187}}
\put(8776,-699){\line( 1,-1){187.500}}
\put(8964,-886){\line( 1, 1){187}}
\put(8964,-1074){\line( 0, 1){188}}
\put(8589,-511){\line( 1,-1){187.500}}
\put(8776,-699){\line( 1, 1){188}}
\put(8589,-136){\line( 1,-1){187.500}}
\put(8776,-324){\line( 1, 1){188}}
\put(8776,-324){\line( 0,-1){375}}
\put(8589,-324){\line( 1,-1){187}}
\put(8776,-511){\line( 1, 1){187.500}}
\put(8401,-1449){\makebox(0,0)[lb]{Suppression of slice number one}}
\put(6901,-1449){\makebox(0,0)[lb]{Staircase }}
\end{picture}
\end{center}
If $E=(E_1,\dots,E_s)$ is a family of staircases, and
$t=(t_1,\dots,t_s)\in \NN^s$,
we put:
\begin{displaymath}
  S(E,t)=(S(E_1,t_1),S(E_2,t_2),\dots,S(E_s,t_s)).
\end{displaymath}
If $(t_1,\dots,t_r)\in (\NN^s)^r$, the recursive formula
\begin{displaymath}
  S(E,t_1,\dots,t_r)=S(S(E,t_1,\dots,t_{r-1}),t_r)
\end{displaymath}
defines the s-tuple of staircases 
$S(E,t_1,\dots,t_r)$ obtained from the s-tuple $E=(E_1,\dots,E_s)$ 
by suppression of $r$ slices in each $E_i$.
\\
If $p\in X$ is a smooth point, a formal neighborhood of $p$
is a morphism $\phi:\s k[[x_1,\dots,x_d]] \fd X$ which induces an 
isomorphism between $\s k[[x_1,\dots,x_d]]$ and the completion
$\widehat {\fxo_p}$ of the local ring of $X$ at $p$. If
$p=(p_1,\dots,p_s)$ is a s-tuple of smooth distinct points, a formal
neighborhood of $p$ is a morphism
$(\phi_1,\dots,\phi_s):U \fd X $ from the disjoint union $U=V_1\coprod
\dots \coprod V_s$ of 
$s$ copies of $\s k[[x_1,\dots,x_d]]$ to $X$, where 
$\phi_i:V_i \fd X$ is a formal neighborhood of
$p_i$. If $D$ is a divisor on $X$, we say that $\phi$ and $D$ are
compatible if $D$ is defined by
the equation $x_1=0$ around each $p_i$ (in particular, $p_i$ is a 
smooth point of $D$).\\
Consider the translation morphism:
\begin{eqnarray*}
Tr_{v_1}: k[[x_1,\dots,x_d]] &\fd & k[[x_1,\dots,x_d]]\ox k[[t]]\\
x_1&\mapsto& x_1\ox 1 -1 \ox t^{v_1}\\
x_i&\mapsto& x_i\ox 1 \mbox{ si $i>1$}
\end{eqnarray*}
If $E_1$ is a staircase, the ideal 
\begin{displaymath}
J(E_1,v_1)=Tr_{v_1}(I^{E_1})k[[x_1,\dots,x_d]] \ox k[[t]]\inc k[[x_1,\dots,x_d]] \ox k[[t]]
\end{displaymath}
defines a flat family $F_1$ of subschemes of
$\s k[[x_1,\dots,x_d]]$ parametrised by $\s k[[t]]$. 
This corresponds geometrically to the family whose fiber over 
$t$ is obtained from $V(I^{E_1})$ by the translation $x_1\mapsto x_1-t^{v_1}$.
If $\phi_1$ is a formal neighborhood of $p_1$, 
$F_1$ can be seen as a flat family of subschemes of $X$ via $\phi_1$,
thus it defines a morphism $\s k[[t]] \fd Hilb(X)$. We denote
by $X(\phi_1,E_1,t,v_1)$ the non closed point of $Hilb(X)$ parametrised by the image of the
generic point. The first coordinate does not play any specific role, 
thus more generally, if $E=(E_1,\dots,E_s)$ is a family of
staircases, if $\phi=(\phi_1,\dots,\phi_s)$ is a formal neighborhood
of $(p_1,\dots,p_s)$, if $v=(v_1,\dots,v_s)\in \NN^s$, one defines
similarly families $F_i\inc X\x \s k[[t]]$ flat over $\s k[[t]]$. 
Since $F_i \cap F_j=\emptyset$ for $i\neq j$, the union 
$F=F_1 \cup \dots \cup F_s$ is still flat over $\s k[[t]]$ and
corresponds to a morphism
$\s k[[t]] \fd Hilb(X)$. We denote by  $X_\phi(E,t,v)$ the image of
the generic point and by $X_\phi(E)=X_\phi(E,0,v)$ the image of
the special point (which does not
depend on $v$).   
Finally, we denote by $[x]$ the integer part of a real $x$. 
\\
We are now ready to state the theorem. By the above,
$\fxl(-X_{\phi}(E,t,v))$ corresponds to a morphism $\s k((t))\fd \GG$ to
a Grassmannian $\GG$, which can be extended to a morphism $\s k[[t]]\fd
\GG$ by valuative properness. The theorem gives a control of the limit
obtained under suitable conditions. 
\begin{thm}
  \label{thm:description de la limite}
  Let $D$ be an effective divisor on a
  quasi-projective scheme $X$,
  $p=(p_1,\dots,p_s)$ be a s-tuple of smooth points of $X$, $\phi$ a formal
  neighborhood of $p$ compatible with $D$,
  $v=(v_1,\dots,v_s)\in \NN^s$ a speed vector, $E=(E_1,\dots,E_s)$ be
  staircases and $X_{\phi}(E,t,v)$ the generic union of subschemes
  defined by $\phi$. Suppose that one can find integers $n_1>\dots>n_r$ 
  such that:
  \begin{listecompacte}
    \item $\forall k$,
  $n_k-n_{k+1}\geq max(v_i)$,
\item $\forall i$, $1\leq i \leq r$,
  $\fxl(-(i-1)D-Z_i)=\fxl(-iD)$
  \end{listecompacte}
where 
  $t_i=([\frac{n_i}{v_1}],\dots,[\frac{n_i}{v_s}] )$,
  $T_i=T(E,t_i)$ and  $Z_i=X_\phi(T_i)$. 
Then 
  \begin{displaymath}
    \lim_{t\fd 0} \fxl(-X_{\phi}(E,t,v)) \subset \fxl(-rD-X_{\phi}(S(E,t_1,\dots,t_r)))
  \end{displaymath}
\end{thm}

\begin{rem}
\label{rem:comparaison alex-hirscho}
  The main lemma 2.3 of
  \cite{alexander-hirschowitz00:grosPointsSuffisammentNombreux}
corresponds essentially to the above theorem with $r=1$. 
\end{rem}
If $X$ is irreducible, $X(E)$ is well defined and it 
specializes to $X_{\phi}(E,t,v)$. Thus we get by
semi-continuity the inequality
\begin{displaymath}
  \dim \fxl(-X(E))\leq \dim \fxl(-X_{\phi}(E,t,v))=\dim \lim_{t\fd 0}\fxl(-X_{\phi}(E,t,v)).
\end{displaymath}
Combining this inequality with the theorem, we obtain the 
following estimate of $ \dim
\fxl(-X(E))$ in terms of a linear system of smaller degree. 
\begin{coro}
  \begin{math}
     \dim \fxl(-X(E)) \leq \dim \fxl(-rD-X_{\phi}(S(E,t_1,\dots,t_r)))
  \end{math}
\end{coro}
\begin{rem}
  In case $\fxl$ is infinite dimensional, the theorem still makes
  sense since Grassmannians of finite codimensional vector spaces of
  $\fxl$ are
  still well defined and the limit makes sense in such a Grassmannian.
\end{rem}

\section{Proof of theorem \ref{thm:description de la limite}}
\label{sec:proof-theorem}
We start with an informal explanation of the ideas in the proof
in the case $s=1$. 
Suppose that we have of family of sections $s(t)$ of $L$ which vanish 
on a moving punctual subscheme $Z(t)=X_{\phi}(E,t,v)$ 
whose support $p(t)$ tends to
$p(0)$ as $t$ tends to $0$. 
Using local coordinates around $p(0)$, the sections of $L$ can 
be considered as functions and the vanishing on $Z(t)$
translates to $s(t)\in J(t)$ where $J(t)$ is the ideal of $Z(t)$. 
Denote by $J_{n_1}$ the restriction of $J(t)$ to the infinitesimal 
neighborhood $\s k[t]/t^{n_1}$ of $t=0$. Suppose that the family 
of sections over  $\s k[t]/t^{n_1}$ is a family of sections 
which vanish on $Z_1$. Then it is a family of sections vanishing on
$D$  since by hypothesis a section
which vanish on $Z_1$ automatically vanishes on $D$. 
If $D$ is defined locally by the equation $x_1=0$, 
this means that $s(t)=x_1s'(t)$ with $s'(t)\in (J_{n_1}:x_1)$. 
Restrict now to the smaller infinetismal neighborhood $\s
k[t]/t^{n_2}$. Suppose that over this restriction, the family 
of sections, which already vanish on $D$, vanish also on 
$Z_2$ (i.e. $s'(t)$ is a family of sections vanishing on $Z_2$). 
Then by hypothesis, the sections vanish twice on $D$. 
Using local coordinates, this means that $s(t)=x_1^2 s''(t)$ 
with $s''(t)\in ((J_{n_1}:x_1)_{n_2}:x_1)$. 
After several restrictions, 
we put $t=0$ and we get $s(0)=x_1^r s^{(r)}(0)$ where $s^{(r)}(0)$ 
is in a prescribed ideal. The control we get in this way of the 
element $s(0)\in \lim_{t\fd 0}\fxl(-X_\phi(E,t,v))$ translates 
into the inclusion  
 \begin{displaymath}
    \lim_{t\fd 0} \fxl(-X_{\phi}(E,t,v)) \subset \fxl(-rD-X_{\phi}(S(E,t_1,\dots,t_r)))
  \end{displaymath}
given by the theorem. 
\\
To play the above game, one needs to be able to compute 
in the successive steps ideals like $((J_{n_1}:x_1)_{n_2}:x_1)$ defined using 
restrictions and transporters.
In view of this explanation, one 
can understand the conditions on the $n_i$ of the theorem 
as follows. The condition $n_1\geq n_2\geq n_3\dots $ 
comes from the  fact that we restrict successivly to smaller 
and smaller neighborhoods. The condition $n_k-n_{k+1}\geq max(v_i)$ 
is a technical condition to be able to compute the successive ideals 
defined via transporters and restrictions. 

Let us start the proof itself now. 
In the context of the theorem, we are given a set of staircases
$E=(E_1,\dots,E_s)$, a vector $v=(v_1,\dots,v_s)$,
a divisor $D$ and a formal neighborhood $\phi$ of $(p_1,\dots,p_s)$ 
in which $D$ 
is given by the equation $x_1=0$ around each $p_i$. For $n>0$, 
we put $R_n=k[[x_1,\dots,x_d]]^s \ox k[[t]]/(t^n)$ and $R_{\infty}=
k[[x_1,\dots,x_d]]^s \ox k[[t]]$. We denote by 
$\psi_{np}:R_n \fd R_p$ the natural projections, which exist for 
$p\leq n \leq \infty$.  If $J\inc R_{\infty}$ is an ideal, we define
we define recursively the ideals  $J_{n_1:n_2:\dots:n_k}\inc R_{n_k}$
and $J_{n_1:n_2:\dots:n_k:}\inc R_{n_k}$ using transporters and 
restrictions by the formulas
\begin{listecompacte}
  \item
  $J_{n_1}=\psi_{\infty n_1}(J)$,
  \item $J_{n_1:n_2:\dots:n_k:}=(J_{n_1:n_2:\dots:n_k}:x_1)$, 
  \item $J_{n_1:n_2:\dots:n_k}=\psi_{n_{k-1}
    n_k}(J_{n_1:n_2:\dots:n_{k-1}:})$
\end{listecompacte}
As explained above, the vector space $\fxl(-X(\phi,E,t,v))$
corresponds 
to a morphism $\s k((t))\fd \GG$ (where $\GG$  is a Grassmannian of
subvectors spaces of $\fxl$) which  extends to a morphism $\s k[[t]]\fd
\GG$. The universal family over the Grassmannian $\GG$ pulls back to a family $U \subset
\s k[[t]]\x \fxl$. Let $e_i$ be a local generator of $L$ at $p_i$. 
Any section $\sigma$ of the line bundle $L$ can be written down 
$\sigma=\sigma_i e_i$ around $p_i$
for some $\sigma_i \in k[[x_1,\dots,x_d]]$. The map:
\begin{eqnarray*}
  \fxl &\fd& k[[x_1,\dots,x_d]]^s\\
  \sigma &\mapsto & (\sigma_1,\dots,\sigma_s)
\end{eqnarray*}
identifies $U$ with a subscheme of $\s k[[t]]\x
k[[x_1,\dots,x_d]]^s$. 
The theorem will be proved if we show that the special fiber $U(0)$
contains only sections vanishing $r$ times on $D$ and if, 
in local coordinates, 
$U(0)$ is included in $x_1^r I^{S(E,t_1,\dots,t_r)}$. 
% (There is a slight abuse of notation here: $x_1^r
% I^{S(E,t_1,\dots,t_r)}\inc k[[t]]\ox k[[x_1,\dots,x_d]]$ 
% whereas $U(0)\inc \s k[[t]]\x k[[x_1,\dots,x_d]]$. 
% In fact, if $I\inc  k[[t]]\ox k[[x_1,\dots,x_d]]$, we still 
% denote by $I$ the subscheme of $ \s k[[t]]\x k[[x_1,\dots,x_d]]$
% whose fiber over a point $p$ is the ideal $I(p)\inc \s k(p)\x k[[x_1,\dots,x_d]]$). 
\\
Let us denote by
$U_{n_i}$ the restriction of $U$ over the subscheme $\s k[[t]]/t^{n_i}$. We show by
induction that:
\begin{displaymath}
  \forall i\geq 1,\ U_{n_i}\inc x_1^i J_{n_1:n_2:\dots:n_i:}
\end{displaymath}
where $J=J(E_1,v_1)\oplus \dots \oplus J(E_s,v_s) \inc R_\infty$.
The fibers of $U$ contain sections of $\fxl$ which vanish on 
$X_{\phi}(E,t,v)$. Since $J$ is the ideal of  $X_{\phi}(E,t,v)$,
this implies the inclusion
$U \inc J$, hence $U_{n_1} \inc J_1$. By
corollary \ref{coro:trace du residuel}, this inclusion implies that 
the fibers of $U_{n_1}$ are
elements of $\fxl$ which vanish on $Z_1$, hence they vanish on $D$ by
hypothesis. It follows that elements of $U_{n_1}$ are dividible by $x_1$
and we can then write: $U_{n_1} \inc x_1 J_{n_1:}$\ . Suppose now that 
 $U_{n_i}\inc x_1^i J_{n_1:n_2:\dots:n_i:}$\ . Then  $U_{n_{i+1}}\inc
 x_1^{i} J_{n_1:n_2:\dots:n_i:n_{i+1}}$. By corollary
 \ref{coro:trace du residuel}, this inclusion implies that the fibers of $U_{n_{i+1}}$ are
elements of $\fxl(-iD)$ which vanish on $Z_{i+1}$, hence they vanish on $D$ by
hypothesis.  It follows that elements of $U_{n_{i+1}}$ are dividible by $x_1^{i+1}$
and we can write $U_{n_{i+1}} \inc x_1^{i+1}
J_{n_1:n_2:\dots:n_{i+1}:} $\ . 
This ends the induction on $i$. In particular, for $i=r$, using
corollary \ref{coro:specialisation du residuel} for the last equality,
we have the required inclusion:
\begin{displaymath}
U(0)=U_{n_r}(0)\inc x_1^r J_{n_1:n_2:\dots:n_r:}(0)= x_1^r I^{S(E,t_1,\dots,t_r)}.\findem
\end{displaymath}
We now turn to the proof of the corollaries \ref{coro:trace du residuel}
and \ref{coro:specialisation du residuel} on which the above proof
relies. Note that 
$J=(J^1,\dots,J^s)$ and 
$I^{T_k}=((I^{T_k})^1,\dots,(I^{T_k})^s)$ are defined componentwise,
the component number $i$ corresponding to the study around the point
$p_i$. 
Thus corollary \ref{coro:trace du residuel} and 
 \ref{coro:specialisation du residuel}
 below can
be proved for each component and one may suppose $s=1$ to prove it. 
We thus suppose for the rest of this section that $s=1$,
that $E=(E_1,\dots,E_s)$ is a staircase given by a height function
$h$, and that $v=(v_1,\dots,v_s)\in \NN$. 
\\
Let $\fxb$ (resp. $\fxc$) be the set of elements $m=(m_2,\dots,m_d)\in
\NN^{d-1}$ such that $h(m)\neq 0$ (resp. $h(m)=0$). Remark that $\fxb$
is finite due to the finitness of $E$.
We denote by 
\begin{listecompacte}
  \item $C(t)\inc R_n$ the $k[[x_1]]\ox k[[t]]$ sub-module containing the
  elements  $\sum a_{m_1m_2\dots m_d}
  x_1^{m_1}x_2^{m_2}\dots x_d^{m_d}\ox f(t)$, where $f(t)\in
  k[[t]]/t^n$ and $(m_2,\dots,m_d)\in \fxc$
\item $C(0)\inc R_1=k[[x_1,\dots,x_d]]$ the   $k[[x_1]]$ sub-module containing the
  series  $\sum a_{m_1m_2\dots m_d}
  x_1^{m_1}x_2^{m_2}\dots x_d^{m_d}$ where $(m_2,\dots,m_d)\in \fxc$
\item $B(m)\inc R_n$ the $k[[x_1]]\ox k[[t]]$ sub-module generated by
  $f_m=(x_1-t^v)^{h(m)}x_2^{m_2}\dots,x_d^{m_d}$, 
\item $B(m,0)\inc R_1=k[[x_1,\dots,x_d]]$ the $k[[x_1]]$ sub-module generated by
  $f_m(0)=(x_1)^{h(m)}x_2^{m_2}\dots,x_d^{m_d}$, 
\item $B_{n_1n_2\dots n_k}(m)\inc R_n$ the $k[[x_1]]\ox k[[t]]$ sub-module generated by
  the elements $f_m,\frac{t^{\alpha_{k-i+1}}f_m}{x_1^i}$, $1\leq i \leq k$, where
 $\alpha_i=max(0,n_i-vh(m))$ for $i>0$. In particular, for $k=0$, 
$B_{n_1n_2\dots n_k}(m)=B(m)$.
\end{listecompacte}
To simplify the notations, we have adopted above the same notation for 
distinct submodules (leaving in distinct ambiant modules).
The following lemma says that the module $B_{n_1n_2\dots n_k}(m)$
is well defined as a sub-module of $R_{j}$ for $j\leq n_{k}$. 
\begin{lm}
  \label{lm:estimation exposant en t des generateurs}
  Let $j\leq n_{k}$. If $i\leq k$, 
  the element $\frac{t^{\alpha_{k-i+1}}f_m}{x_1^i} \in R_j$.
  In particular $B_{n_1n_2\dots n_k}(m) \inc R_j$ is
  well defined for  $j\leq n_{k}$. If in addition, $j \leq n_{k+1}$,
  then  $\frac{t^{\alpha_{k-i+1}}f_m}{x_1^i}$ is a multiple of $x_1$.
\end{lm}
\demo. First, if $l<i$, the coefficient of $x_1^l$ in $t^{\alpha_{k-i+1}}
  f_m$ is a multiple of $t^{\alpha_{k-i+1}}t^{v(h(m)-l)}$. This term 
is zero in $R_j$ since the exponent of $t$ is at least $n_{k-i+1}-vl
  \geq n_{k}+(i-1)v-vl
  \geq n_{k}\geq j$. It follows that  $\frac{t^{\alpha_{k-i+1}}
  f_m}{x_1^i}\in R_j$ is well defined. A similar estimate shows that 
for $l\leq i$,  the coefficient of $x_1^l$ in $t^{\alpha_{k-i+1}}
  f_m$ is zero in $R_j$ for $j \leq n_{k+1}$. Thus
  $\frac{t^{\alpha_{k-i+1}}f_m}{x_1^i}$
is a multiple of $x_1$.
\findem

\begin{lm} \label{description de I^E gradue}
  \begin{listecompacte}
    \item As $k[[x_1]]$-modules, $I^E=\bigoplus_{m \in \fxb}B(m,0)
    \oplus C(0) \inc k[[x_1,\dots,x_d]]$
  \item As $k[[x_1]]\ox k[[t]]$-modules, $J=\bigoplus_{m \in \fxb}B(m)
    \oplus C(t) \inc R_n$
  \end{listecompacte}
\end{lm}
\demo: This is a straightforward verification left to the reader. 
\findem

\begin{lm}  \label{lm:description des residuels}
  We have the equality of $k[[x_1]]\ox k[[t]]$-modules:
  \begin{listecompacte}
    \item $J_{n_1:\dots :n_k}=\bigoplus _{m \in \fxb} B_{n_1n_2\dots
    n_{k-1}}(m) \oplus C(t)\inc R_{n_k}$
\item $J_{n_1:\dots :n_k:}=\bigoplus _{m \in \fxb} B_{n_1n_2\dots
    n_k}(m) \oplus C(t) \inc R_{n_k}$
  \end{listecompacte}
\end{lm}
\demo. Let us say that the number of indexes of $J_{n_1:\dots :n_k}$
and $J_{n_1:\dots :n_k:}$ is respectivly $2k-1$ and $2k$. We prove
the lemma by induction on the number $i$ of indexes. If $i=1$,
we get from the preceding lemma the equality
\begin{eqnarray*}
  J_{n_1}=\psi_{\infty n_1}(J)&=&\sum_{m\in \fxb} \psi_{\infty n_1}(B(m))
  +  \psi_{\infty n_1}(C(t))\\
  &=& \sum_{m\in \fxb}B(m) + C(t)\mbox { in } R_{n_1}.
\end{eqnarray*}
The last sum is obviously direct, thus it is the required equality. \\
Suppose now that we want to prove the lemma for $i=2k-1$. This is
exactly the same reasoning as in the case $i=1$, substituting 
$ J_{n_1:\dots :n_k}$, $J_{n_1:\dots :n_{k-1}:}$ and
$\psi_{n_{k-1}n_k}$ for $J_{n_1}$, 
$J$, and $\psi_{\infty,n_1}$. \\
For the last case $i=2k$. Taking the transporter from the expression
of $ J_{n_1:\dots :n_k}$ coming from induction hypothesis, we get:
\begin{displaymath}
J_{n_1:\dots :n_k:}=\bigoplus _{m \in \fxb} (B_{n_1n_2\dots
    n_{k-1}}(m):x_1) \oplus (C(t):x_1)
\end{displaymath}
The equality $(C(t):x_1)=C(t)$ is obvious, so we are done if we prove 
the equality  $(B_{n_1n_2\dots
    n_{k-1}}(m):x_1)= B_{n_1n_2\dots
    n_{k}}(m)$ in the ambiant module $R_{n_k}.$
The inclusion $\supset$ is clear since for every generator $g$ of $ B_{n_1n_2\dots
    n_{k}}(m)$, $x_1g$ is a multiple of one of the generators of $B_{n_1n_2\dots
    n_{k-1}}(m)$. As for the reverse inclusion, if $z \in (B_{n_1n_2\dots
    n_{k-1}}(m):x_1)$, one can write down 
  \begin{displaymath}
    x_1z=\sum_{1 \leq i \leq k-1}P_i
    \frac{t^{\alpha_{k-i+1}}f_m}{x_1^i} +x_1P_0f_m+Q_0f_m\ \  (*)
  \end{displaymath}
where $P_i\in k[[x_1]]\ox k[[t]]$ and $Q_0 \in k[[t]]$. By lemma 
\ref{lm:estimation exposant en t des generateurs},
the terms $ \frac{t^{\alpha_{k-i+1}}f_m}{x_1^i} \in R_{n_k}$ are dividible by
$x_1$, thus $x_1$ divides $Q_0f_m$.  
It follows  that the coefficient $Q_0t^{vh(m)}x_2^{m_2}\dots
x_d^{m_d}$ of $x_1^0$ in $Q_0 f_m$ is zero, which happens only if $Q_0$ is 
a multiple of $t^{max(0,n_k-vh(m))}=t^{\alpha_{k}}$. Writing down 
$Q_0=t^{\alpha_{k-1+1}}$ and 
dividing the
displayed equality $(*)$ by $x_1$ shows that $z\in  B_{n_1n_2\dots
    n_{k}}$, as expected. \findem

\begin{coro}
  \label{coro:trace du residuel}
   $J_{n_1:n_2:\dots:n_k}\inc I^{T_k}$ 
\end{coro}
\demo. In view of the previous lemma, and since the inclusion $C\inc
I^{T_k}$ is obvious, 
one simply has to check that 
the generators of $B_{n_1:n_2:\dots:n_k}(m)$ 
verify the inclusion. The generators are explicitly given thus
this is a straightforward verification.
\findem

\begin{coro}
  \label{coro:specialisation du residuel}
   $J_{n_1:n_2:\dots:n_k}(0)=I^{S(E,t_1,\dots,t_k)}$.
\end{coro}
\demo. According to lemmas \ref{lm:description des residuels} 
and \ref{description de I^E gradue}, it suffices to show that
$B_{n_1n_2\dots
    n_{k}}(m,0) \inc k[[x_1]]$ is the submodule generated by 
$x_1^{h(m)- p(m)}$ where $p(m)$ is the number of $t_i$'s verifying
$t_i<h(m)$. Since the generators of $B_{n_1n_2\dots
    n_{k}}(m)$ are explicitely given, 
the corollary just comes from the evaluation of these
 generators at $t=0$. 
\findem

\section{The Hilbert function of $k^2$ fat points in $\plp$}
\label{sec:hilbert-function}

In this section, we compute the Hilbert function of the generic union 
of $k^2$ fat points in $\plp$ of the same multiplicity $m$. \\
We work over a field of characteristic 0. 
\begin{defi}
  If $Z\inc \plp$ is a zero-dimensional subscheme of degree $deg(Z)$, we denote
  by $H_v(Z):\NN\fd \NN$ the virtual Hilbert function of $Z$ defined by the
  formula $H_v(Z,d)=min(\frac{(d+1)(d+2)}{2},deg(Z) )$. The critical
  degree for $Z$, denoted by $d_c(Z)$ is the smallest integer $d$ such
  that $H_v(Z,d)>deg(Z)$. 
\end{defi}

\begin{thm} \label{thm:fonction de Hilbert de k^2 points}
  Let $Z$ be the generic union of $k^2$ fat points of multiplicity
  $m$. Then $H(Z)=H_v(Z)$. 
\end{thm}

Let us recall the following well known lemma:

\begin{lm}
  \label{lm:les degres critiques suffisent}
  If $H(Z,d) \geq H_v(Z,d)$  for $d= d_c(Z)$ and $d=d_c(Z)-1$, then $H(Z)=H_v(Z)$. 
\end{lm}

\begin{defi}
The regular staircase $R_m\inc \NN^2$ is the set defined by the
relation $(x,y)\in R_m \Leftrightarrow x+y<m$. A quasi-regular
staircase $E$ is a staircase such that $R_m\subset E \subset R_{m+1}$
for some $m$. A right specialized staircase is a staircase such that
$((x,y)\in E \ and\ y>0) \Rightarrow (x+1,y-1)\in E$. A monomial
subscheme of $\plp$ with staircase $E$ is a punctual subscheme
supported by a point $p$ which is defined by the ideal $I^E$ in 
some formal neighborhood of $p$. 
\end{defi}

Our first intermediate goal is lemma 
\ref{lm:elimination des points simples}
which says that under suitable conditions,
if $Z=L\cup R\inc \plp$ is a subscheme with $L$ included in a line, 
the Hilbert function of $Z$ is determined by that of $R$.  

\begin{prop}
  Let $Z$ be a generic union of fat points. The following
  conditions are equivalent.
  \begin{listecompacte}
    \item $H(Z)=H_v(Z)$
    \item  there exists a quasi-regular right-specialized staircase $E$ and a collision
  $C$ of the fat points which is monomial with staircase $E$.
\item  there exists a quasi-regular staircase $E$ and a collision
  $C$ of the fat points which is monomial with staircase $E$.
  \end{listecompacte}
\end{prop}
\demo. 
$1\Rightarrow 2$. 
Let $\rho_t$ be the automorphism of $\plp=Proj(k[X,Y,H])$
defined 
for $t \neq 0$ by 
$f_t:X\mapsto \frac{X}{t} , Y\mapsto \frac{Y}{t},H \mapsto H$.
Consider the collision $C=\lim_{t\fd 0} f_t(Z)$. It is a
subscheme of the affine plane $\s k[x=\frac{X}{H} ,y=\frac{Y}{H} ]$
supported by the origin $(0,0)$. 
It is shown in 
\cite{evain97:lesCollisionsDeterminentLaPostulation-CRAS} that
if $H(Z)=H_v(Z)$, then there is an integer $m$ such that the ideal of 
$C$ verifies $I^{R_{m+1}}\inc I(C) \inc I^{R_m}$. Thus
$I(C)=V\oplus k[x,y]_{\geq m+1}$ where $ k[x,y]_{\geq m+1}$
stands for the vector space generated by the monomials of 
degree at least $m+1$, and $V\inc k[x,y]_{m}$. 
Let now $g_t:x\mapsto x-ty , y\mapsto y$. Then the ideal of 
$D=\lim_{t \fd \infty}g_t(C)$ is $I(D)=W\oplus k[x,y]_{\geq m+1}$
where $W=\lim_{t \fd \infty}g_t(V)$ is a vector space which admits 
a base of the form $y^m,xy^{m-1},\dots,x^ky^{m-k}$. Thus $I(D)=I^E$
for some quasi-regular right-specialized staircase $E$. And $D$ is a
collision of the fat points since it is a specialisation of the
collision $C$ and since being a collision is a closed condition. 
\\
$2\Rightarrow 3$ is obvious. \\
$3\Rightarrow 1$. 
If there exists a collision $C$ associated with a
quasi-regular staircase $E$, then by semi-continuity $H(Z,d)\geq
H(C,d)=min(\frac{(d+1)(d+2)}{2}, \#E )= min(\frac{(d+1)(d+2)}{2},deg(C)
)= min(\frac{(d+1)(d+2)}{2},deg(Z))=H_v(Z,d)$. Since the 
well known reverse inequality
$H_v(Z,d) \geq H(Z,d)$ is always true, we have the required equality
$H_v(Z,d) = H(Z,d)$. \findem

\begin{lm}
  \label{lm:elimination des points simples}
  Let $R\inc \plp$ be a generic union of fat points, $D\inc plp$ be a generic
  line, $L\inc D$ be a subscheme whose support  is generic in
  $D$. Let $Z=R \cup L$ and suppose that the degree of $L$ satisfies
  $deg(L)\leq d_c(R)$. Then $H(R)=H_v(R)$ implies
  $H(Z)=H_v(Z)$. 
\end{lm}
\demo. By the above lemma and its proof, there exists a quasi-regular
right specialized staircase $E$ and a collision $C$ 
of the fat points 
supported by the origin of $\AA^2=\s k[x,y]$ such that 
the ideal of $C\inc \pla$ is
$I(C)=I^E$. By the genericity hypothesis, $L$ can be specialized to 
the subscheme $L(t)$ with equation $(y-t,x^{deg(L)})$. Obviously 
$L(t)$ is monomial with staircase $F=\{(0,0),(1,0),\dots,(r,0)\}$.
Let $D=\lim_{t\fd 0}C\cup L(t)$. 
By
\cite{hirscho85:methodeHoraceManuscripta},
$I(D)=I^G$ for some monomial staircase $G$.
Moreover, the explicit description of $G$ given
in \cite{hirscho85:methodeHoraceManuscripta}
( $G$ is the ``vertical collision'' of $E$ and $F$)
shows that $G$ is quasi-regular.
Since $Z=R\cup L$ can be specialized to a scheme $D$ defined by 
a quasi regular staircase, $H(Z)=H_v(Z)$.
\findem

\begin{lm}
  \label{lm:estimation du degre critique}
  Let $Z\inc \plp$ be a union of $k^2$ fat points of multiplicity $m$ with $k
  \geq 4$. The
  critical degree $d_c(Z)$ verifies $km+1<d_c(Z)\leq km+k-2$.
\end{lm}
\demo: Direct calculation. \findem

\noindent
\textit{Proof of theorem \ref{thm:fonction de Hilbert de k^2 points}.}\\
We show by induction on $k$ that the Hilbert function of the 
generic union $Z$ of $k^2$ fat
points of multiplicity $m$ is the virtual Hilbert function $H_v(Z)$.
If $k\leq 3$, this is known by
\cite{nagata60:9points_imposent_conditions_indep}.
So we may suppose $k\geq 4$.  According to lemma \ref{lm:les degres
  critiques suffisent}, we only need to check that $H(Z,d) \geq H_v(Z,d)$
for $d=d_c(Z)$ or $d=d_c(Z)-1$, and, by lemma \ref{lm:estimation du
  degre critique}, such a $d$ verifies $d=km+s$ for some $s$ satisfying $0\leq s\leq
k-2$.
By semi-continuity, 
it suffices to specialize $Z$ to a scheme $Z'$ with
$H(Z',d)\geq H_v(Z,d)$. First, we choose 
a generic line $D$ and generic points $p_1,\dots,p_{2k-1}$ on $D$.
We divide the $k^2$ fat points into three subsets $E_1,E_2,E_3$ of
respective cardinal
$k,k-1,(k-1)^2$. 
We specialize the $k$ fat points of $E_1$
on the points $p_{k},\dots,p_{2k-1}$. We leave the generic $(k-1)^2+(k-1)$ 
points of $E_3\cup E_2$ in their generic position.
We denote by $\fxl$ the set of sections of $\fxo(d)$ which vanish on the
fat points of $E_1 \cup E_3$. Since the points of $E_1$ have been
specialised, we have  by semi-continuity the
inequality:
\begin{displaymath}
 (*)\ H(Z,d) \geq \frac{(d+1)(d+2)}{2} -\dim \fxl(-X(E))
\end{displaymath}
where 
\begin{displaymath}
E=(\underbrace{R_m,\dots,R_m}_{(k-1)\ copies}).
\end{displaymath}
We now make a further specialisation, moving the $k-1$ 
fat points of $E_2$ on the points
$p_1,\dots,p_{k-1}$ using theorem \ref{thm:description de la limite}.
To this end, we fix the notations. We choose a formal neighborhood $\phi$
of  $p=(p_1,\dots,p_{k-1})$, a number $N>>0$ and we take the
speed vector
\begin{displaymath}
  v=(\underbrace{N,\dots,N}_{k-s-2\ times},\underbrace{N+1,\dots,N+1}_{s+1\ times}).
\end{displaymath}
Finally, we let 
\begin{displaymath}
n_i=(N+1)(m-i+1)-1, 1\leq i\leq m.
\end{displaymath}
Let us check that the conditions of 
theorem \ref{thm:description de la limite} apply.
The condition $n_k-n_{k+1}\geq max(v_i)$ is obviously satisfied.
As for the remaining condition,
remark that $\fxl(-(i-1)D)$ is a set of sections
of $\fxo(d-i+1)$ which vanish on $p_{k}^{m-i+1},\dots,p_{2k-1}^{m-i+1}$.
In particular, if $Z_i$ is a punctual subscheme of
$D$ of cardinal $d-i+2-k(m-i+1)=s+1+(i-1)(k-1)$ whose support does not meet 
the union $p_k \cup \dots \cup p_{2k-1}$, then $\fxl(-iD-Z_i)=\fxl(-(i+1)D)$.
In our case, 
$Z_i$ is a union of one-dimensional fat points of the line $D$. Let us compute its  degree.  
The subscheme $Z_i$ is supported by $p_1 \cup \dots \cup
p_{k-1}$ and we denote by $d_j$ the degree of the part $(Z_i)_{p_j}$ supported by $p_j$.
It is the cardinal $m-[\frac{n_i}{v_j}]$ 
of the slice $T(R_m,[\frac{n_i}{v_j}])$, that is $d_j=i-1$
if $j\leq k-s-2$ and $d_j=i$ if $k-s-1\leq j \leq k-1$. Thus the
degree of $Z_i$ is the sum of the $d_j$, that is  $s+1+(i-1)(k-1)$.
We can then apply theorem  \ref{thm:description de la limite}
and its corollary.
We conclude that:
\begin{displaymath}
(**)\ \dim \fxl(-X(E)) \leq \dim
\fxl(-mD-X_{\phi}(S(E,t_1,\dots,t_m)).
\end{displaymath}
The linear system $\fxl(-mD)$ is the set of sections of
$\fxo(d-m)$ which vanish on the union $Z'$ of the fat points of
$E_3$. Moreover, 
$X_{\phi}(S(E,t_1,\dots,t_m))$ is the union $L$ of the 
one-dimensional fat points
of $p_1^m\cap D,\dots,p_{k-s-2}^m \cap D$.  It follows that 
\begin{displaymath}
   (***)\ \dim
\fxl(-mD-X_{\phi}(S(E,t_1,\dots,t_m))=\frac{(d-m)(d-m+1)}{2} H(Z'\cup L,d-m)
\end{displaymath}.
By lemma \ref{lm:elimination des points simples} and the induction, we
have 
\begin{displaymath}
  (****)\ H(Z'\cup L,d-m)=H_v(Z'\cup L,d-m)
\end{displaymath}
Now, by construction (or by an easy direct calculation), 
\begin{displaymath}
  (*****)\ H_v(Z'\cup L,d-m)-\frac{(d-m)(d-m+1)}{2}
  =H_v(Z,d)-\frac{d(d+1)}{2} 
\end{displaymath}
Putting together the displayed equalities and inequalities (*)\dots
(*****) gives the
required inequality $H(Z,d)\geq H_v(Z,d)$. 
\findem

\section{Collisions of fat points}
\label{sec:collision-fat-points}

We start with a definition of a generic successive collision
of fat points in $\pla$. We proceed 
by induction. A generic successive collision of 
one fat point $p^m$ is the fat point itself. Suppose defined the 
generic successive collision $Z_{m_1\dots m_{k-1}}$
of  $p_1^{m_1},\dots ,p_{k-1}^{m_{k-1}}$. 
Let $C(d)$ be the generic curve of degree $d$ containing the 
support $O$ of $Z_{m_1\dots m_{k-1}}$. Let 
\begin{displaymath}
Z_{m_1\dots
    m_{k}}(d)=\lim_{p\in C(d),\ p\fd O} Z_{m_1\dots m_{k-1}}\cup p^{m_k}.
\end{displaymath}

\begin{prop}\label{prop:defDeLaCollisionSuccessiveGenerique}
  There exists an integer $d_0$ such that $\forall d\geq d_0$,
  $Z_{m_1\dots m_{k}}(d)=Z_{m_1\dots m_{k}}(d_0)$. We denote this 
  subscheme by $Z_{m_1\dots m_{k}}$ and this is by definition the
  generic successive collision of $p_1^{m_1},\dots,p_k^{m_k}$. 
\end{prop}
\demo. Consider the morphism $f:\pla\setminus\{O\}\fd Hilb(\pla)\inc Hilb(\plp)$ which
sends the point $p\in \pla$ to the subscheme $Z_{m_1\dots m_{k-1}}
\cup p^{m_k}$. It extends to a morphism $\tilde f: S\fd Hilb(\plp)$, where
$\pi:S\fd \pla$ is a composition of blowups (of simple points). 
The embeddings  $\s k[t]/(t^d) \fd \pla$ sending the support of 
$\s k[t]/(t^d)$ to $O\in \pla$ form an irreducible variety and we
denote by $g:\s k[t]/t^d \fd \pla$ the corresponding generic embedding.
For $p\geq d$, the intersection $C(p)\cap O^d$ of the curve 
with the fat point is isomorphic as an
abstract scheme to $\s k[t]/(t^d)$; since for any embedding 
$i:\s k[t]/(t^d) \fd \pla$, there exists a curve of degree $d$ 
which contains the image $Im(i)$, it follows  that $C(p)\cap O^d$ is 
the subscheme associated with the generic embedding $g$.
In particular, $C(p)\cap O^{d_0}=C(d_0)\cap O^{d_0}$ if $p\geq d_0$. 
Choose $d_0>n$ where $n$ is the number of blowups in $\pi$. Since
the order of contact of $C(p)$ and $C(d_0)$ is at least $d_0$, the
number of blowups is not sufficient to separate the curves and the
strict transforms $\tilde C(p)\inc S$ and $\tilde C(d_0)\inc S$ intersect in a
point $s$. It follows that  $Z_{m_1\dots m_{k}}(p)=Z_{m_1\dots
  m_{k}}(d_0)=\tilde f(s)$.
\findem
Our goal is to compute the generic collision $Z_{mmmm}$ of $4$ fat
points of multiplicity $m$. 
\begin{rem}
  With the notations of proposition
  \ref{prop:defDeLaCollisionSuccessiveGenerique}, the integers $d_0$
  which appear in the definition of $Z_{mmmm}$ will always be equal to
  $1$. In other words, the collision will be shown to 
  depend only on the tangent
  directions of the approaching fat points. 
\end{rem}
We will describe $Z_{mmmm}$ as a pushforward via a blowup 
$\pi:\tilde S \fd \pla$, where $\pi$ is 
the blowup defined by  the 
following Enriques diagram
  
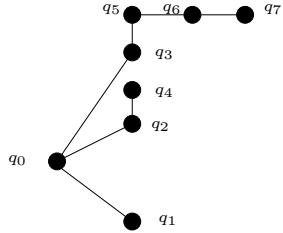
\begin{figure}[h] 
     \begin{center}
        \input{constellation.pstex_t}
     \end{center} 
  \end{figure}.\\ 
We recall for convenience what this means. 
Let $q_0\in\pla$,
$q_1,q_2,q_3$ be three distinct tangent directions at $q_0$.
Let 
\begin{displaymath}
  \eta:S_1\fd S_0=\pla
\end{displaymath}
be the blowup of $q_0$, and $Q_0\inc S_1$ the exceptional divisor.
Let 
\begin{displaymath}
  S_2\fd S_1
\end{displaymath}
be the blowup of $(q_1\cup q_2 \cup q_3) \inc Q_0$, and $Q_1,Q_2,Q_3
\inc S_2$ the respective exceptional divisors. 
If $Q_i\inc S_{n_i}$ is an exceptional divisor, and if $S_j \fd
S_{n_i}$ is a sequence of blowups, we still denote by $Q_i\inc S_j$ (resp. we
denote by $E_i\inc S_j$) the strict transform (resp. the total transform) of
$Q_i$ in $S_{j}$. With this convention, let $q_4=Q_0 \cap Q_2 \in S_2$, $q_5 = Q_0 \cap
Q_3 \in S_2$. Let 
\begin{displaymath}
S_3 \fd S_2
\end{displaymath}
be the blowup of $q_4 \cup q_5$,
$Q_4,Q_5$ the corresponding exceptionnal divisors. 
Let 
 $q_6=Q_3\cap Q_5 \in S_3$, 
$S_4\fd S_3$ the blowup of $q_6$, $Q_6$ its exceptional divisor.
Let  $q_7=Q_6 \cap
Q_3 \in S_4$ and $\tilde S=S_5\fd S_4$ the blowup 
of $q_7$. 
We denote by 
\begin{displaymath}
  \rho:\tilde S \fd S_1\ \ and\ \ \pi:\tilde S \fd \pla
\end{displaymath}
the compositions of  the blowups introduced above. 
As explained, each point $q_i$ defines a divisor 
$E_i \inc \tilde S$. 
If $(m_0,\dots,m_7)\in \NN^8$, the ideal $\pi_*(\fxo_{\tilde S}(-\sum m_iE_i))$
is a punctual subscheme supported by $q_0$ which we will represent graphically 
with a label $m_i$ at the point of the Enriques diagram corresponding
to $q_i$. For instance, the subscheme  $\pi_*(\fxo_{\tilde
  S}(-8E_0-2E_1-E_2-E_4-3E_3))$
is associated  with the following diagram.

  
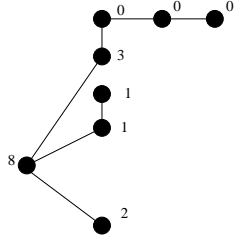
\begin{figure}[h] 
     \begin{center}
        \input{exempleSchemaComplet.pstex_t}
     \end{center} 
  \end{figure}\\
The following theorem describes the successive 
collision of four fat points which approach on curves  $C_i$ with distinct 
tangent directions. This includes in particular the generic successive
collision.
\begin{thm}
  Let $q_0 \in \pla$, $q_1,q_2,q_3$ three distinct tangent directions at
  $q_0$ and $C_1,C_2,C_3$ be three smooth curves passing through $p_0$
  with tangent direction $q_1,q_2,q_3$.
  Let 
  $Z_{mmmm}$ be the collision of the fat points
  $p_0^m,p_1^m,p_2^m,p_3^m$ where:
  \begin{listecompacte}
    \item $p_0$ is located at $q_0$,
    \item $p_1$ moves on the curve $C_1$ (resp. $p_2$ on $C_2$, $p_3$
    on $C_3$).
  \end{listecompacte}
Then $Z_{mmmm}$ 
  is defined by the following Enriques diagram, which depends on $m$
  modulo 4.
  
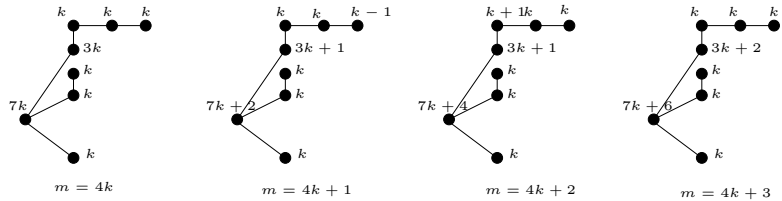
\begin{figure}[h] 
     \begin{center}
        \input{4collision.pstex_t}
     \end{center} 
  \end{figure}
\end{thm}
\demo. All cases are similar and we prove the theorem in the case
$m=4k$. We choose a formal neighborhood $\xi$ of $p=(q_1,q_2,q_3)\in
(S_1)^3$ such that $Q_0\inc S_1$ is defined by the equation $x_1=0$
around each $q_i$ and such that
$C_3$ is defined by $x_2=0$ around $q_3$ (this is possible since $C_3$
is smooth). 
Let $n=(m-1,m-5,\dots,3)$. Let $F_m$ be the staircase defined by
the height function $h_{F_m}(d)=h_{R_m}([\frac{d}{2}])$, and let 
$G=S(R_m,n)$ be the staircase obtained from $R_m$ by suppression of
the slices indexed by $n$. Let $X_{\xi}(R_k,F_k,G_k)\inc S_1$ be the
subscheme defined by the formal neighborhood $\xi$ and the staircases
$R_k,F_k,G_k$. 
According to the correspondance between complete ideals and monomial
subschemes formulated in \cite{evain99:4hgrosPoints}, if the $m_i$'s
are the integers defined in the Enriques diagram,
\begin{displaymath}
  \rho_* \fxo_{\tilde S}(-\sum m_i E_i))=
  \fxo_{S_1}(-m_0Q_0-X_{\xi}(R_k,F_k,G))) \ \ (*)
\end{displaymath}
Let $J(p_3)$ denote the ideal of $Z_{mmm}\cup p_3^m$. I claim that 
we are done if we prove the inclusion 
\begin{displaymath}
  \lim_{p_3 \fd p_0} \eta^* J(p_3)\inc
  H^0(\fxo_{S_1}(-m_0Q_0-X_{\xi}(R_k,F_k,G))\ \ (**).
\end{displaymath}
Indeed, we would then have the inclusions
\begin{eqnarray*}
  I_{Z_{mmmm}}&\inc & \eta_* \eta^*  I_{Z_{mmmm}}=
 \eta_* \eta^*  \lim_{p_3 \fd p_0} J(p_3)\\
  &\inc &  \eta_*  \lim_{p_3 \fd p_0} \eta^* J(p_3)\\
  &\inc &  \eta_*  H^0(\fxo_{S_1}(-7kQ_0-X_{\phi}(R_k,F_k,G))\ \ by\ (**)\\
  &\inc &  H^0(\eta_*  (\fxo_{S_1}(-7kQ_0-X_{\phi}(R_k,F_k,G))\\
  &\inc &  H^0(\eta_*  \rho_*  \fxo_{\tilde S}(-\sum m_i E_i))\ \ by\ (*) \\
  &\inc & I_Z \mbox{ where $I_Z=\pi_*\fxo_{\tilde S}(-\sum m_iE_i)$}.
\end{eqnarray*}
According to 
\cite{casas90:pointsInfinimentVoisinsMathAnnalen}, 
since the Enriques diagram defining $Z$ 
is unloaded, $deg(Z)=\sum \frac{m_i(m_i+1)}{2}$
which is immediatly checked to be $4
\frac{4k(4k+1)}{2}=deg(Z_{mmmm})$. Summing up, $Z$ and $Z_{mmmm}$ are two
punctual subschemes of the same degree with $ I_{Z_{mmmm}} \inc I_Z$,
thus they are equal.
\\
It remains to prove the displayed inclusion $(**)$ using our theorem. By
\cite{evain_these} or \cite{walter95:collisionsDeTroisGrosPointsUnpublished},
\begin{displaymath}
  \eta^*I_{Z_{mmm}}=H^0 \fxo_{S_1}(-6kQ_0-X_{\psi}(R_{2k},F_{2k}))
\end{displaymath}
where $\psi$ is the formal neighborhood of $(q_1,q_2)$ induced by the
formal neighborhood $\xi$ of $(q_1,q_2,q_3)$. Thus
\begin{displaymath}
  \lim_{p_3\fd p_0}\eta^*J(p_3)=\lim_{t\fd 0} \fxl(-X_{\phi}(R_m,t,v=1))
\end{displaymath}
where $\phi$  is the formal neighborhood of $p_3$ induced by the
formal neighborhood $\xi$ of $(q_1,q_2,q_3)$ and
$\fxl=H^0(\fxo_{S_1}(-6kQ_0 -X_{\psi}(R_{2k},F_{2k})))$. 
To apply theorem 
\ref{thm:description de la limite} with $X=S_1$, $s=1$,$D=Q_0$, and
$n=(m-1,m-5,\dots,3)$, the verification $\fxl((-i+1)D-Z_i)=\fxl(-iD)$
is needed. Elements of $\fxl((-i+1)D-Z_i)$ are sections of $\fxo_{S_1}((-6k-i+1)Q_0)$
that vanish on 
\begin{displaymath}
X_{\psi}(R_{2k-i+1},F_{2k-i+1})\cup
Z_i=X_{\xi}(R_{2k-i+1},F_{2k-i+1},T(R_m,m-1-4(i-1))).
\end{displaymath}
Since the
intersection 
\begin{displaymath}
Q_0 \cap
X_{\xi}(R_{2k-i+1},F_{2k-i+1},T(R_m,m-1-4(i-1)))
\end{displaymath}
has
degree $3(2k-i+1)+(4i-3)$  greater than the degre $6k+i-1$ of the restriction
$\fxo_S((-6k-i+1)Q_0)_{|Q_0}$, it follows that any section of
$\fxl((-i+1)D-Z_i)$ vanishes on $D$. Thus we can apply the theorem 
and we get:
\begin{eqnarray*}
  \lim_{t \fd 0}\fxl(-X_{\phi}(R_m,t,1))&\inc&
\fxl(-kQ_0-X_{\phi}(S(R_m,n)))\\
&\ &\hspace{2cm} =\\
&& H^0(\fxo_{S_1}(-7kQ_0-X_{\psi}(R_k,F_k)-X_{\phi}(S(R_m,n))))\\
&&\hspace{2cm} =\\
&& H^0(\fxo_{S_1}(-m_0Q_0-X_{\xi}(R_k,F_k,S(R_m,n))),
\end{eqnarray*}
which concludes the proof. \findem

\bibliography{/users/evain/perso/fichiersConfig/modeles/texBiblio.bib}
\bibliographystyle{plain} 

\end{document}